\documentclass[12pt]{amsart}

\newcounter{EGnum}
\newcommand{\EGNumber}{\theEGnum\stepcounter{EGnum}}

\newenvironment{EG}[1]{{\vspace{1 ex}}\noindent {\sc Example \EGNumber.}{#1}{\hfill{$\diamondsuit$}}\\{\vspace{1 ex}}}

\usepackage{latexsym} 
\usepackage{amsmath} 
\usepackage{amssymb} 
\usepackage{amscd} 
\usepackage{psfig} 
\usepackage{multicol}
\usepackage{amsfonts}

\newtheorem{theorem}{Theorem}

\newtheorem{thm}{Theorem}[section] 
\newtheorem{corollary}[thm]{Corollary} 
\newtheorem{lemma}[thm]{Lemma} 
\newtheorem{prop}[thm]{Proposition} 
\newtheorem{definition}[thm]{Definition} 
\newtheorem{remark}[thm]{Remark}

\newcommand{\cpth}{\ensuremath{\mathbb{C}P^3}}
\newcommand{\cptw}{\ensuremath{\mathbb{C}P^2}}

\def\to{\rightarrow} 
\def\into{\hookrightarrow}

\newcommand{\C}{\mathbb C}
\newcommand{\algg}{\ensuremath{\mathfrak{g}}}

\newcommand{\algl}{\ensuremath{\mathfrak{l}}}
\newcommand{\algh}{\ensuremath{\mathfrak{h}}}
\newcommand{\algs}{\ensuremath{\mathfrak{s}}}

\def\iso{\cong}

\def\g{\mathfrak{g}}

\title[Equivariant Cohomology and Residual Actions]{The equivariant cohomology of Hamiltonian $G$-spaces From Residual $S^1$ Actions}
\author{Rebecca Goldin}
\address[Rebecca Goldin]{Department of Mathematics\\University of Maryland\\College Park, MD 20742-4015}
\email{goldin@math.umd.edu}
\author{Tara S. Holm}
\address[Tara Holm]{Department of Mathematics\\Massachusetts Institute of Technology 2-251\\Cambridge MA 02139}
\email{tsh@math.mit.edu}

\begin{document}

\begin{abstract}
We show that for a Hamiltonian action of a compact torus $G$ on a
compact, connected symplectic manifold $M$, the $G$-equivariant
cohomology is determined by the residual $S^1$ action on the
submanifolds of $M$ fixed by codimension-1 tori.  This theorem allows
us to compute the equivariant cohomology of certain manifolds, which
have pieces that are four-dimensional or smaller.  We give several
examples of the computations that this allows.
\end{abstract}

\maketitle

\section{Introduction}\label{se:intro}

It has long been a ``folk theorem'' that, for a Hamiltonian torus
action on a symplectic manifold, the associated equivariant cohomology
is determined by $S^1$ actions on certain submanifolds. Recently,
Tolman and Weitsman \cite{TW:CohomHamTspaces} used equivariant Morse
theory to prove that the cohomology is determined by that of the {\em
one-skeleton}, the subspace given by the closure of all points whose
orbit under the torus action is one-dimensional. Here we use a
powerful result of Chang and Skjelbred
\cite{CS} to give a short proof
of a slightly stronger statement. 

Let $M$ be a compact, connected symplectic manifold, and let $G$ be a
compact torus acting on $M$ in a Hamiltonian fashion. Let 
$M^G$ denote the fixed point set of the action. The inclusion
$r:M^G\hookrightarrow M$ induces a map 
$$
r^*: H_G^*(M)\longrightarrow H_G^*(M^G).
$$ 
F. Kirwan \cite{Kirwan:cohom} proved that $r^*$ is an
injection. We find a simple description of the image of $r^*$ in terms
of the $S^1$-equivariant cohomology of submanifolds of $M$ fixed by
any codimension-1 subtorus of $G$.

We first fix some notation. Let $H\subset G$ be a codimension-1 torus
and let $M^H$ denote its fixed point set. Let $X_H$ denote a connected
component of $M^H$.  Let 
$$
r_{X_H}^*:H_G^*(X_H)\rightarrow H_G^*(X_H^G) 
$$ 
be the map induced by the inclusion of the fixed
points $X_H^G$ into $X_H$ and 
$$ 
i^*_{X_H}: H_G^*(M^G)\rightarrow H_G^*(X_H^G) 
$$ 
be the map induced by the inclusion of $X_H^G =
X_H\cap M^G$ into $M^G$.  We can reduce the computation of
$G$-equivariant cohomology of $M$ to the computation of
$S^1$-equivariant cohomology of submanifolds of $M$ as follows.

\begin{theorem}\label{th:fromCS}
Let $M$ be a compact, connected symplectic manifold with a Hamiltonian
action of a torus $G$.
Let $r^*:H_G^*(M)\longrightarrow H_G^*(M^G)$ be the map induced by the
inclusion of the fixed point set. A class $f\in H_G^*(M^G)$ is in the
image of $r^*$ if and only if
$$
i_{X_H}^*(f)\in r_{X_H}^*(H_G^*(X_H)).
$$
for all codimension-1 subtori $H\subset G$ and connected components
$X_H$ of $M^H$.
\end{theorem}

In particular, consider the case in which $G$ acts with isolated fixed
points, and $\dim X_H\leq 2$ for all $X_H$.  Theorem \ref{th:fromCS}
gives an explicit description of $r^*$, which was proven in
significant generality in \cite{GKM:eqcohom}. See, for
example, \cite{GS:supersymmetry}, Chapter 11.  First we find the
cohomology of each component $X_H$.  If $\dim X_H =2$, then $X_H$ is
diffeomorphic to $S^2$ with a Hamiltonian $S^1\iso G/H$ action with
fixed points denoted $\{N,S\}$. Suppose first that $G\iso S^1$ and
$H=\{ 0\}$. In that case,
$$ 
r_{X_H}^*:H^*_{S^1}(S^2)\hookrightarrow H_{S^1}^*(\{N,S\}) 
$$ 
is the inclusion induced by $\{N,S\}\subset S^2$. By the Atiyah-Bott
Berline-Vergne localization theorem \cite{AB:local}, \cite{BV:local},
any element $f\in
H_{S^1}^*(\{N,S\})$ in the image of $r_{X_H}^*$ must satisfy the
property that
\begin{equation}\label{eq:diff}
f_N-f_S \in x\cdot \C[x]
\end{equation}
where $f_N$ and $f_S$ are the restrictions of $f$ to the points $N$
and $S$, respectively. We have identified the equivariant
cohomology of a point $H_{S^1}^*(pt)$ with $\C[x]$.  

 Let $R$ be the graded
ring $H_{S^1}^*(\{N,S\})$ subject to the above restriction. A quick
dimension check shows that as modules over $H_{S^1}^*(pt)$,
$H_{S^1}^*(S^2) = R$. However, the module structure forces the rings
to be equal, so that condition (\ref{eq:diff}) is the only condition
for $f\in im(r_{X_H}^*)$.

Suppose now that $G\iso T^n$.

\begin{prop}
Suppose that $S^2$ is a Hamitonian $G$-space for $G\iso T^n$.
Let $H$ be a codimension $1$ subtorus which acts trivially.  Then a
function $f=(f_N,f_S)\in S(\algg^*)\oplus S(\algg^*)$ is in the image
of $r^*:H_G^*(S^2)\to H_G^*(\{ N,S\} )$ if and only if
$$
f_N-f_S\in\ker (\pi_H),
$$
where $\pi_H:S(\algg^*)\to S(\algh^*)$ is induced by the projection
$\algg^*\to\algh^*$. 
\end{prop}

\begin{proof}
Choose a complement $L$ to $H$ in $G$ and write $S(\algl^*)\iso
\C[x]$.  We note that $H_G^*(S^2)=H_L^*(S^2)\otimes H_H^*(pt)$ and
$H_G^*(\{ N,S\} )=H_L^*(\{ N,s\} )\otimes H_H^*(pt)$ because $H$ acts
trivially on $S^2$.  Furthermore, the map
$$
r^*:H_L^*(S^2)\otimes S(\algh^*)\to H_L^*(\{ N,s\} )\otimes S(\algh^*),
$$
is the identity of the second component, where $S(\algh^*)$ is the
$H$-equivariant cohomology of a point.  Thus, by (\ref{eq:diff}),
$f\in H_G^*(\{ N,S\} )$ is in the image of 
$r^*H_G^*(S^2)\to H_G^*(\{ N,S\} )$ if and only if
$f_N-f_S\in(x)\cdot\C[x]\otimes S(\algh^*)$.  But this is precisely
the kernel of $\pi_H$.
\end{proof}

Using this description of $H_{G}^*(X_H)$ in the case that $\dim
(X_H)\leq 2$,, we have the following corollary, a theorem of Goresky,
Kottwitz and MacPherson \cite{GKM:eqcohom}.

\begin{corollary}\label{cor:GKM}
\cite{GKM:eqcohom}
Let $M$ be a compact, symplectic manifold with a Hamiltonian action of
a compact torus $G$. Asume that $M^G$ consists of isolated fixed
points $\{p_1,\dots,p_d\}$ and that each component $X_H$ of $M^H$ has
dimension 0 or 2 for $H\subset G$ a codimension-1 torus. Let $f_{i}$
be the restriction of a class $f\in H_G^*(M)$ to the fixed point
$p_i$. Let $\pi_H:\mathfrak{g}^*\rightarrow \mathfrak{h}^*$ be the
projection induced by the inclusion $\algh\into\algg$. Then the map 
$$
r^*:H_G^*(M)\longrightarrow H_G^*(M^G) =\bigoplus_{p\in M^G}H_G^*(pt) 
$$
has image $(f_{1},\dots, f_{d})$ such that $$
\pi_H(f_{i})=\pi_H(f_j)
$$
whenever $\{p_i,p_j\}=X_H\cap M^G$.

\end{corollary}

This theorem can be stated in terms of graphs (cf. \cite{GZ:graphs}).
The natural generalization of graphs are hypergraphs. These are
explored by the authors in \cite{GH:hypergr}.

We now discuss the more general case, which extends the above
corollary. We allow the fixed point sets of codimension $1$
subtori of $G$ to have dimension 0, 2, or 4.  The equivariant
cohomology of $M$ in this case can be computed as follows.

\begin{theorem}
\label{th:4total}
Suppose that $M$ is a
compact, connected symplectic manifold with an effective Hamiltonian
$G$ action.  Suppose further that the $G$ action has only isolated
fixed points $M^G=\{ p_1,\dots ,p_d\}$ and that $\dim X_H\leq 4$ for
all $H\subset G$ of codimension $1$ and $X_H$ a connected component of $M^H$.
As before, let $f_i\in H_G^*(pt)$ denote the restriction of $f\in
H_G^*(M)$ to the fixed point $p_i$.
The image of the injection $r^*:H_{G}^*(M)\rightarrow H_{G}^*(M^G)$ is 
the subalgebra of functions $(f_1,\dots
,f_d)\in\bigoplus_{i=1}^d  S(\g^*)$ which satisfy
$$
\left\{\begin{array}{ll}
\pi_H(f_{i_j})=\pi_H(f_{i_k}) & \mbox {if } \{ p_{i_1},\dots
,p_{i_l}\}=X_H^G  \\
\sum_{j=1}^l\frac{f_{i_j}}{\alpha_1^{i_j}\alpha_2^{i_j}}\in
S(\g^*) & \mbox {if }\{ p_{i_1},\dots ,p_{i_l}\}=X_H^G \mbox{ and }
\dim X_H =4
\end{array}
\right . 
$$
for all $H\subset G$ codimension-1 tori, where $\alpha_1^{i_j}$ and $\alpha_2^{i_j}$ are the (linearly dependent) weights of the $G$ action on $T_{p_{i_j}}X_H$.
\end{theorem}

The results we present here are not generally true outside the
symplectic setting.  We rely heavily on the fact that the restriction
map to the equivariant cohomology of the fixed point set is an
injection.  We also make use of the Morse theory associated to a Morse 
function obtained from the moment map for the Hamiltonian action.

In Section \ref{se:maintheoremproof} we prove
Theorem~\ref{th:fromCS}. In Section~\ref{se:dim4} we prove
Theorem~\ref{th:4total}.  Lastly, in Section~\ref{se:example} we give 
several examples of the computations allowed by this theorem.

\section{Reduction to the study of circle actions}\label{se:maintheoremproof}

The proof of Theorem~\ref{th:fromCS} is based on the following result of
Chang and Skjelbred \cite{CS}.
The description given here is due to Brion and Vergne
\cite{BrV:local}.

\begin{lemma}\cite{CS}\label{le:cs}
The image of $r^*:H_G^*(M)\rightarrow H_G^*(M^G)$ is the set

$$
\bigcap_Hr_{M^H}^*(H_G^*(M^H)),
$$

\noindent where the intersection in $H_G^*(M^G)$ is taken over all
codimension-one subtori $H$ of $G$ and $r_{M^H}$ is the inclusion of
$M^G$ into $M^H$. 
\end{lemma}

\begin{remark}
\label{rm:finite}
In fact, the only nontrivial contributions to this intersection are those
codimension-one subtori $H$ which appear as isotropy groups of elements of
$M$.  Since $M$ is compact, there are only finitely many such isotropy
groups. 
\end{remark}

We will now restate and prove the theorem reducing the computation of
$H_G^*(M)$ to computations of $H_{S^1}^*(X)$ for various submanifolds
$X\subseteq M$.

\medskip

\noindent {\bf Theorem~\ref{th:fromCS}.}{\em \ \ 
Let $r^*:H_G^*(M)\longrightarrow H_G^*(M^G)$ be the map induced by the
inclusion of the fixed point set. Let $i_{X_H}^*:H_G^*(M^G)\to
H_G^*(X_H^G)$ and $r_{X_H}^*:H_G^*(X_H)\to H_G^*(X_H^G)$ be the maps
defined in Section~\ref{se:intro}. A class $f\in H_G^*(M^G)$ is in the
image of $r^*$ if and only if
$$
i_{X_H}^*(f)\in r_{X_H}^*(H_G^*(X_H)).
$$
for all codimension-1 subtori $H\subset G$ and connected components
$X_H$ of $M^H$.
}

\begin{proof}
By Lemma~\ref{le:cs}, $f\in im(r^*)$ if and only if $f$ is in the
intersection of $r_{M^H}^*(H_G^*(M^H))$ over all codimension-one subtori
$H$, where $r_{M^H}^*:H_G^*(M^H)\to H_G^*(M^G)$.  Equivalently,
$$ 
f\in \bigcap_{H}r^*_{M^{H}}(\bigoplus_{X_H}H_G^*(X_H)),
$$
where the direct sum is taken over all connected components $X_H$ of
$M^H$. Let $k_{X_H}:H_G^*(X_H)\rightarrow H_G^*(M^H)$ be the map which
extends
any class on $X_H$ to 0 on other components of $M^H$. Let
$k_{X_H^G}:H_G^*(X_H^G)\rightarrow H_G^*(M^G)$ be the same map on the
fixed point sets. Then
$$
r^*_{M^{H}}(\bigoplus_{X_H}H_G^*(X_H))=\bigoplus_{X_H}r_{M^H}^*\circ
k_{X_H}(H_G^*(X_H)).
$$
As $k_{X^G_H}\circ r^*_{X_H}=r^*_{M^H}\circ k_{X_H}$, we have that
$f$ is in $im(r^*)$ if and only if
\begin{equation}
\label{eq:main1}
f\in \bigoplus_{X_H}k_{X_H^G}\circ r_{X_H}^*(H_G^*(X_H)),
\end{equation}
for all $H$.  
Now note that $i^*_{X_H}\circ k_{X^G_H}=id$.  Because the $X_H$ are
disjoint, we can now apply $i_{X_H}^*$ to (\ref{eq:main1}) to
get
\begin{equation}
\label{eq:main2}
i_{X_H}^*(f)\in r_{X_H}^*(H_G^*(X_H)),
\end{equation}
for every $H$ and $X_H$.  However, since $\bigoplus_{X_H}i_{X_H}^*$ is
an injection, we can apply $\bigoplus_{X_H}k_{X_H^G}$ to
(\ref{eq:main2}) to get (\ref{eq:main1}). Thus, (\ref{eq:main2}) and
(\ref{eq:main1}) are equivalent.  This completes the proof.
\end{proof}
 
This provides another proof for a result of Tolman and Weitsman
\cite{TW:CohomHamTspaces}. 
\begin{definition}
Let $N\subset M$ be the set of points whose orbits under the $G$ action are 1-dimensional. The {\em one-skeleton} of $M$ is the closure $\overline{N}$.
\end{definition}
Tolman and Weitsman show that the image of $r^*:H_G^*(M)\rightarrow
H_G^*(M^G)$ is equal to the image of the cohomology of the
one-skeleton.
\begin{thm}\cite{TW:CohomHamTspaces}
Let $M$ be a compact symplectic manifold with a Hamiltonian torus
action by $G$. Let $M^G$ be the fixed point set, and $\overline{N}$ be
the one-skeleton. Let $r:M^G\hookrightarrow M$ be the inclusion of the
fixed point set to $M$ and $j: M^G \hookrightarrow \overline{N}$ be
the inclusion to $\overline{N}$.
The induced maps $r^*:H_G^*(M)\rightarrow
H_G^*(M^G)$ and $j^*:H_G^*(\overline{N})\rightarrow H_G^*(M^G)$ on
equivariant cohomology have the same image.
\end{thm}

\begin{proof}
Because $G$ acts effectively, $N$ consists of points fixed by some
codimension-1 torus $H\subset G$ but not by all of $G$, i.e.
$$
N = \bigcup_H M^H \backslash M^G
$$
where the union is taken over all codimension-1 tori $H\subset
G$. As noted above, this is a finite union over
all codimension-1 $H$ which appear as isotropy subgroups of points in
$M$.
Then $\overline{N}= \bigcup_H M^H$,  and the inclusion
$\gamma_H: M^H\hookrightarrow M$ factors through the inclusion
$\gamma: \overline{N}\hookrightarrow M$ for each codimension-1 torus
$H$ in $G$. It follows that the induced maps in cohomology also
factor. Furthermore, there is an inclusion 
$$
H_G^*(\overline{N})\hookrightarrow \bigoplus_{i=1}^k H_G^*(M^{H_i}),
$$
where $H_i, i=1,\dots,k$ are the codimension-1 tori which appear as
isotropy subgroups of $G$.  Theorem \ref{th:fromCS} implies that the
map $r^*:H_G^*(M)\rightarrow H_G^*(M^G)$ factors through the map
$$
\bigoplus_{i=1}^k r_{M^{H_i}}^*:\bigoplus_{i=1}^k H_G^*(M^{H_i})\longrightarrow
H_G^*(M^G)
$$
But then $r^*$ must factor through $j^*:H_G^*(\overline{N})\rightarrow
H_G^*(M^G)$.
\end{proof}

\vspace{.1in}

Now suppose that $M^G$ consists of isolated fixed points.
Then $$H_G^*(M^G)= \bigoplus_{p\in M^G}S(\mathfrak{g}^*)$$ and any $f\in
H_G^*(M^G)$ is a map
$f:M^G\rightarrow S(\mathfrak{g}^*)$. Furthermore, as $X_H$ and
$X_H^G$ have trivial $H$ actions, we can rewrite Theorem
\ref{th:fromCS} in the following way. 

\begin{thm}\label{th:fromCSisolatedfp}
Under the above hypotheses, the image of $r^*$ is the set of
$f:M^G\rightarrow S(\mathfrak{g}^*)$ such that
$$
i_{X_H}^*(f):X_H^G\rightarrow
S(\mathfrak{g}^*)$$ 
 is in the image of 
$$
r_{X_H}^*:H_{G}^*(X_H)\rightarrow H_{G}^*(X_H^G)= \bigoplus_{p\in
X_H^G}S(\mathfrak{g}^*),
$$
where $r_{X_H}^*$ is restriction for each fixed point. 
\end{thm}

\section{An extension of a theorem of GKM}\label{se:dim4}

We use Theorem~\ref{th:fromCSisolatedfp} to compute the equivariant
cohomology of $M$ in the case in which $\dim X_H\leq 4$ for all
codimension-1 tori $H\subset G$ and connected components $X_H$ of the
fixed point set $M^H$. In Section~\ref{se:intro} we considered the
case that $\dim X_H\leq 2$.

First, let $X$ be a compact, connected symplectic
four-manifold with an effective Hamiltonian $G=S^1$ action with
isolated fixed points.  Then the equivariant
cohomology can be computed as follows.

\begin{prop}
\label{pr:GKMin4} Let $X$ be a compact, connected symplectic
4-manifold with an effective Hamiltonian $S^1$ action with isolated
fixed points $X^{S^1}=\{ p_1,\dots ,p_d\}$.
The map $r^*:H_{S^1}^*(X)\rightarrow H_{S^1}^*(X^{S^1})$ induced by
inclusion is an injection with image
\begin{equation}
\label{eq:S^1onM^4}
\{ (f_1,\dots ,f_d)\in\bigoplus_{i=1}^d S(\algs^*)\ |\ f_i-f_j\in
x\cdot\C [x],\ \sum_{i=1}^d\frac{f_i}{\alpha_1^i\alpha_2^i}\in
S(\algs^*)\},
\end{equation}
where $\alpha_1^i$ and $\alpha_2^i$ are the (linearly
dependent) weights of the $S=S^1$ isotropy action on $T_{p_i}X$.
\end{prop}

\begin{proof}
The map $r^*$ is injective because $X$ is equivariantly formal.  We
know that the $f_i$ must satisfy the first condition because the
functions constant on all the vertices are the only equivariant
classes in degree $0$, as $\dim H_{S^1}^0(X)=1$.  The second
condition is necessary as a direct result of the Atiyah-Bott
Berline-Vergne (ABBV) localization theorem
(\cite{AB:local},\cite{BV:local}).   Notice that this condition
gives us one relation in degree $2$ cohomology.  A dimension count
shows us that these conditions are sufficient.  As an $S(\algs^*)$-module, $H_{S}^*(X)\iso
H^*(X)\otimes H_{S}^*(pt)$.  Thus, the equivariant
Poincar\'{e} polynomial is
\begin{eqnarray*}
P^S_t(X) & = & (1+(d-2)t^2+t^4)\cdot (1+t^2+t^4+\dots) \\
         & = & 1+(d-1)t^2+dt^4+\dots+dt^{2n}+\cdots.
\end{eqnarray*}

\noindent  As $H_{S}^*(X)$ is generated in degree 2, the $d-1$
degree 2 classes given by the $(f_1,\dots,f_d)$ subject to the ABBV
condition generate the
entire cohomology ring. Thus we have found all the conditions.
\end{proof}

We now prove a slightly more general proposition. 

\begin{prop}
\label{pr:GKMgenin4} Let $X$ be a compact, connected  symplectic
4-manifold with a Hamiltonian $G$
action with isolated fixed points $X^G=\{ p_1,\dots,p_d\}$.   Suppose
further that there is a 
codimension-$1$ subtorus $H$ which acts trivially.
The map $r^*:H_{G}^*(X)\rightarrow H_{G}^*(X^{G})$ induced by
inclusion is an injection with image
\begin{equation}
\label{eq:T^nonM^4}
\{ (f_1,\dots ,f_d)\in\bigoplus_{i=1}^d S(\algg^*)\ |\ f_i-f_j\in
\ker (\pi_H),\ \sum_{i=1}^d\frac{f_i}{\alpha_1^i\alpha_2^i}\in
S(\algg^*)\},
\end{equation}
where $\pi_H$ is the map $S(\algg^*)\to S(\algh^*)$, and $\alpha_1^i$
and $\alpha_2^i$ are the (linearly dependent) weights of the $G$
isotropy action on $T_{p_i}X$. 
\end{prop}

\begin{proof}
As in the case where $X\iso S^2$,
$$
H_G^*(X)=H_{G/H}^*(X)\otimes S(\algh^*).
$$
Furthermore,
$$
H_G^*(X^G)=H_{G/H}^*(X^G)\otimes S(\algh^*).
$$
Again, choose a
complement $L$ to $H$, and write $S(\algl^*)\iso
\C[x]$. Then $H_{G/H}^*(X^G)$ can be identified with $\bigoplus_{p\in
X^G}\C[x]$.  
By Proposition~\ref{pr:GKMin4}, we have $f\in H_G^*(X^G)$ is in the image of
$r^*:H_G^*(X)\to H_G^*(X^G)$ if and only if
the component of $f$ in $H_{G/H}^*(X^G)\iso\bigoplus_{p\in X^G}\C[x]$
satisfies the conditions (\ref{eq:S^1onM^4}) of
Proposition~\ref{pr:GKMin4}.   But then $f$ must
satisfy the conditions (\ref{eq:T^nonM^4}).
\end{proof}

We now discuss the more general case, which extends the result due to
\cite{GKM:eqcohom} (Corollary \ref{cor:GKM}).  Suppose that $M$ is a
compact, connected symplectic manifold with an effective Hamiltonian
$G$ action.  Suppose further that this $G$ action has only isolated
fixed points $M^G=\{ p_1,\dots ,p_d\}$ and that $\dim X_H\leq 4$ for
all $H\subset G$ and $X_H$ a connected component of $M^H$, as above.
As before, let $f_i\in H_G^*(pt)$ denote the restriction of $f\in
H_G^*(M)$ to the fixed point $p_i$.  The equivariant cohomology of $M$
can be computed as follows.

\vskip 0.1in

\noindent {\bf Theorem~\ref{th:4total}.}{\em \ \ 
Under the above hypotheses, the image of the injection $r^*:H_{G}^*(M)\rightarrow H_{G}^*(M^G)$ is the subalgebra of functions $(f_{1},\dots ,f_{d})\in\bigoplus_{i=1}^d S(\g^*)$ which satisfy
$$
\left\{\begin{array}{ll}
\pi_H(f_{i_j})=\pi_H(f_{i_k}) & \mbox {if } \{ p_{i_1},\dots
,p_{i_l}\}=X_H^G  \\
\sum_{j=1}^l\frac{f_{i_j}}{\alpha_1^{i_j}\alpha_2^{i_j}}\in
S(\g^*) & \mbox {if }\{ p_{i_1},\dots ,p_{i_l}\}=X_H^G \mbox{ and }
\dim X_H =4
\end{array}
\right . 
$$
for all $H\subset G$ codimension-1 tori, where $\alpha_1^{i_j}$ and $\alpha_2^{i_j}$ are the (linearly dependent) weights of the $G$ action on $T_{p_{i_j}}X_H$.
}

\begin{proof}
By Theorem~\ref{th:fromCSisolatedfp}, $im(r^*)$ consists of
$(f_1,\dots, f_d)$ which have certain properties restricted to each
$X_H$. Proposition~\ref{pr:GKMgenin4} lists these restrictions for
each $X_H$ of dimension $4$. The conditions for $X_H$ of dimension $2$
are discussed in Section~\ref{se:intro}. A quick check shows that
these are exactly the conditions listed above.
\end{proof}

\section{Examples}\label{se:example}
Here we demonstrate the use of Theorem~\ref{th:4total} in computing
equivariant cohomology.
In the first example, we compute the $S^1$-equivariant cohomology of 
$\cptw$ with a Hamiltonian circle action.  In the second, we calculate 
the $T^2$-equivariant cohomology of $\cpth$.  Finally, we find the
$S^1$-equivariant cohomology of a $4$-dimensional manifold obtained by 
symplectic reduction.

\begin{EG}{  
Consider $\C P^2$ with homogeneous coordinates
$[z_0:z_1:z_2]$.  Let $T=S^1$ act on $\cptw$ by 
$$
e^{i\theta}\cdot [z_0:z_1:z_2] = [e^{-i\theta}z_0:z_1:e^{i\theta}z_2].
$$
This action has three fixed points: $[1:0:0]$, $[0:1:0]$ and
$[0:0:1]$.  

The weights at these fixed points are
$$
\begin{array}{cc}
\mbox{Fixed point} & \mbox{Weights} \\
\mbox{$p_1=[1:0:0]$} & x, 2x, \\
\mbox{$p_2=[0:1:0]$} & -x, x, \\
\mbox{$p_3=[0:0:1]$} & -2x, -x,
\end{array}
 $$ where we have identified $\mathfrak{t}^*$ with degree one
polynomials in $\C [x]$. As cohomology elements, these are assigned
degree two. 
The image of the equivariant cohomology $H_{S^1}^*(\cptw
)$ in $H_{S^1}^*(\{ p_1,p_2,p_3\} )\iso \bigoplus_{i=1}^3\C [x]$ is
the subalgebra generated by the triples of functions $(f_1,f_2,f_3)$ such that
$$
f_i-f_j\in x\cdot\C [x]  \mbox{ for every $i$ and $j$, and } 
$$
$$
\frac{f_1}{2x^2}-\frac{f_2}{x^2}+\frac{f_3}{2x^2}\in\C [x]. 
$$}
\end{EG}

We use the cohomology computed above to compute the $T^2$-equivarant cohomology of $\cpth$.

\begin{EG}{  
The second example we consider is a $T^2$ action on $\C P^3$. Consider
$\C P^3$ with homogeneous coordinates $[z_0:z_1:z_2:z_3]$.  Let $T^2$
act on $\cpth$ by 
$$ 
(e^{i\theta_1},e^{i\theta_2})\cdot
[z_0:z_1:z_2:z_3] =
[e^{-i\theta_1}z_0:z_1:e^{i\theta_1}z_2:e^{i\theta_2}z_3].  
$$ 
This
action has four fixed points, $[1:0:0:0]$, $[0:1:0:0]$, $[0:0:1:0]$
and $[0:0:0:1]$.  The image of the moment map for this action is show
in the figure below.

\begin{figure}[h]
\centerline{
\psfig{figure=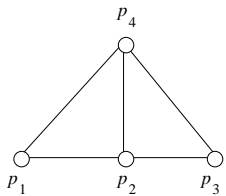,height=1in}
}
\smallskip
\centerline{
\parbox{4.5in}{\caption[Moment map]{{\small This shows the image of 
the moment map for $T^2$ acting on $M=\cpth$, as described above.}}}
}
\end{figure}

The weights at these fixed points are
$$
\begin{array}{cc}
\mbox{Fixed point} & \mbox{Weights} \\
\mbox{$p_1=[1:0:0:0]$} & x, 2x, x+y, \\
\mbox{$p_2=[0:1:0:0]$} & -x, x, y, \\
\mbox{$p_3=[0:0:1:0]$} & -2x, -x, y-x, \\ 
\mbox{$p_4=[0:0:0:1]$} & -x-y, -y, x-y.
\end{array}
$$
Theorem~\ref{th:4total} tells us that the image of the equivariant
cohomology $H_{T^2}^*(\cpth )$ in $H_{T^2}^*(\{ p_1,p_2,p_3,p_4\}
)\iso \bigoplus_{i=1}^4\C [x,y]$ is the ring of functions
$(f_1,f_2,f_3,f_4)$ such that 
$$
f_i-f_j\in (x)\cdot\C [x,y]  \ \  \mbox{ for every $i,j\in \{ 1,2,3\}$,}
$$
\begin{eqnarray*}
\frac{f_1}{2x^2}-\frac{f_2}{x^2}+\frac{f_3}{2x^2}& \in& \C [x,y], \\
f_1-f_4 & \in & (y+x)\cdot\C [x,y], \\
f_2-f_4 & \in & (y)\cdot\C [x,y], \\
f_3-f_4 & \in & (y-x)\cdot\C [x,y].
\end{eqnarray*}

\begin{figure}[h]
\centerline{
\psfig{figure=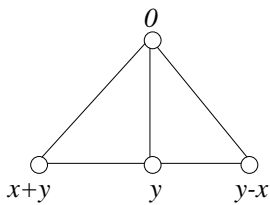,height=1in}
}
\smallskip
\centerline{
\parbox{4.5in}{\caption[Cohomology class]{{\small This shows an equivariant
class of $H_{T^2}^*(\cpth )$, shown as an element of the equivariant
cohomology of the fixed points.}}}
}
\end{figure}}
\end{EG}

\begin{EG}{
Let $\mathcal{O}_\lambda$ be the coadjoint orbit of $SU(3)$ through
the generic point $\lambda\in \mathfrak{t}^*$, the dual of the Lie algebra
$\mathfrak{t}$ of the maximal 2-torus $T$ in $SU(3)$. Recall that $T$
acts on $\mathcal{O}_\lambda$ in a Hamiltonian fashion, and (one
choice of) the moment
map 
$$
\Phi_T: \mathcal{O}_\lambda\longrightarrow \mathfrak{t}^*
$$
takes each matrix to its diagonal entries. Equivalently, $\Phi_T$ is the
composition of the inclusion of $\mathcal{O}_\lambda$ into $\mathfrak{su}(3)^*$ and
projection of $\mathfrak{su}(3)^*$ onto $\mathfrak{t}^*$.

We compute the equivariant cohomology of $M
=\mathcal{O}_\lambda/\!/H$, the {\em symplectic reduction} of
$\mathcal{O}_\lambda$ by a circle $H$ chosen such that the reduced
space is a manifold.  Let $H\subset T$ be any copy of $S^1$ which fixes
a two-sphere in $\mathcal{O}_\lambda$. Then the moment map
$\Phi_H:\mathcal{O}_\lambda\rightarrow \mathfrak{h}^*$ for
the $H$ action is the map $\Phi_T$ followed by the projection
$\pi_H:\mathfrak{t}^*\rightarrow \mathfrak{h}^*$ induced by the
inclusion $\mathfrak{h}\hookrightarrow \mathfrak{t}$. The symplectic
reduction at $\mu$ by $H$ is by definition
$$
M =\mathcal{O}_\lambda/\!/H := \Phi_H^{-1}(\mu)/H,
$$
where $\mu$ is a regular value for $\Phi_H$. Note that there is a
residual $T/H\iso S^1$ action on $M$. We use Theorem \ref{th:fromCS}
to calculate the the corresponding equivariant cohomology of $M$.

\begin{figure}[h]
\label{fig:su3/s^1}
\centerline{
\psfig{figure=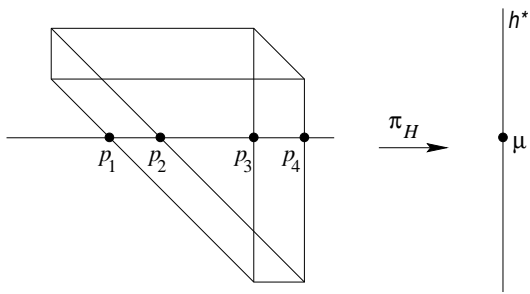,height=1.5in}
}
\smallskip
\centerline{
\parbox{4.5in}{\caption[Moment map]{{\small On the left is the image of 
the moment map for $T$ acting on $\mathcal{O}_{\lambda}$.  The cut
through the moment polytope for $\mathcal{O}_{\lambda}$ corresponds to
the symplectic reduction of $\mathcal{O}_{\lambda}$ by $H$ at $\mu$,
for some choice of $\algh^{\perp}$. 
}}}  }
\end{figure}

One can easily see that there are four fixed points of this action,
which we denote by $p_i$ for $i=1,\dots, 4$. For each $p_i$,
$\Phi_T^{-1}(p_i)$ lies on a two-sphere in $\mathcal{O}_{\lambda}$,
denoted $S^2_i$, which is fixed by a a subgroup $H_i\iso S^1$ of
$T$. Note that $H_i$ is complementary to $H$ in $T$.

The weights of the $T/H$ action on the tangent space $T_{p_i}M$ are
determined by the $T$ action on $S^2_i$.  Let $n_i$ and $s_i$ be fixed
points of the $T$ action on $S^2_i$. Note that the condition that
$\mu$ be a regular value of $\Phi_H$ ensures that $\Phi_T^{-1}(p_i)\neq
n_i, s_i$. Furthermore, by assumption the set $\Phi_T^{-1}(p_i)$ is
point-wise fixed by $H_i$. Thus in the reduction, the $T/H$ action on
$T_{p_i}M$ is isomorphic to the $H_i$ action on this space.

Denote the weights of the $T$ action on $T_{n_i}\mathcal{O}_\lambda$ by
$\pm\alpha_1, \pm\alpha_2$, and $\pm \alpha_3=\pm(\alpha_1+\alpha_2)$,
where the signs depend on $i$. The weights of the $H$ action on the
reduction $M$ are determined
by projecting the $\alpha_i$ to $\mathfrak{h_i}^*$.

At $p_1$, the weight $\alpha_3$ projects
to $0$ and the other two weights both project to
the generator $x$ of $S(\mathfrak{h_1}^*)\iso \C[x]$. Similarly, at $p_2$ the weights are
$x$ and $-x$, at $p_3$ they are $x$ and $-x$ and at $p_4$ they are
both $-x$.  The image of the moment map $\Phi_H: M\rightarrow
\mathfrak{h}^*$, with weights, is shown in
Figure~\ref{fig:redsu3}.

\begin{figure}[h]
\label{fig:redsu3}
\centerline{
\psfig{figure=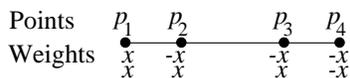,height=0.375in}
}
\smallskip
\centerline{
\parbox{4.5in}{\caption[Moment map]{{\small 
The image of the moment map for the $T/H$ action on
$M=\mathcal{O}_{\lambda}/\!/H$, with the weights for the isotropy
action on the tangent space of the fixed points.
}}}  }
\end{figure}

Finally, this tells us that the equivariant cohomology of $M$ is
\begin{eqnarray*}
H_{S^1}^*(M) & \iso & \{ f:V\rightarrow \C [x]\ |\ f_i-f_j\in
x\cdot\C[x], \\
      & &  \frac{f_1}{x^2} - \frac{f_2}{x^2} -
\frac{f_3}{x^2} + \frac{f_4}{x^2}\in \C[x]\}.
\end{eqnarray*}
Notice that this computation leads us to the
$T/S^1$-equivariant cohomology of $M\iso \mathcal{O}_{\lambda}/\!/S^1$
for a coadjoint orbit of $SU(n)$, as the submanifolds that appear are
identical to those shown the above $SU(3)$ case.
}
\end{EG}

\end{document}